\documentclass[letterpaper,12pt]{article}
\title{\bf Tropical Mathematics}
 \author{David Speyer \  and \
  Bernd Sturmfels\footnote{Supported as a Clay Mathematics Institute  Senior Scholar
           at the IAS/Park City Mathematics Institute, July 2004}
 \\ {\small Department of Mathematics, University of
California,
Berkeley} }

\date{}

\usepackage{amsthm, amsmath, amssymb, amsfonts, graphics}

\theoremstyle{plain}
\newtheorem{thm}{Theorem}

\newtheorem{fact}[thm]{Fact}

\theoremstyle{definition}

\newcommand{\rr}{\mathbb{R}}

\begin{document}
\maketitle

These are the notes for the Clay Mathematics Institute
Senior Scholar Lecture which was
delivered by Bernd Sturmfels in Park City, Utah,
on July 22, 2004. The topic of this lecture is the ``tropical
approach'' in
mathematics, which has gotten a lot of attention recently
in  combinatorics, algebraic geometry and related fields.
It offers an an elementary introduction to this subject,
touching upon Arithmetic, Polynomials, Curves,
Phylogenetics and Linear Spaces. Each section ends with a suggestion
for further research.   The bibliography contains numerous
references for further reading in this field.

The adjective ``tropical''  was coined by French mathematicians,
including Jean-Eric~Pin \cite{Pin}, in the honor of their Brazilian colleague
Imre Simon \cite{Simon}, who was one of the pioneers in
min-plus algebra. There is no deeper meaning in
the adjective ``tropical''. It simply stands for the French view of Brazil.

\section{Arithmetic}

Our basic object of study is the {\em tropical semiring} $\,(\rr \cup
\{\infty\}, \oplus, \odot)$.
As a set this is just the real numbers $\rr$, together with an extra
element $\infty$
which represents infinity. However, we redefine the basic arithmetic
operations of addition and multiplication of real numbers as follows:
$$ x \,\oplus \, y \,\,\,:= \,\,\, {\rm min}(x,y) \qquad \hbox{and} \qquad
x \, \odot \,y \,\,\,:= \,\,\, x + y . $$
In words, the {\em tropical sum} of two numbers is their minimum,
and the {\em tropical product} of two numbers is their sum.
Here are some examples of how to do arithmetic in this
strange number system. The tropical sum of $3$ and $7$ is $3$.
The tropical product of $3$ and $7$ equals $10$. We write this as follows:
$$ 3 \, \oplus \, 7 \,\,\, = \,\,\, 3 \qquad \hbox{and} \qquad
     3 \,\odot \, 7 \,\,\, = \,\,\, 10. $$
 Many of the familar axioms of arithmetic remain valid
 in tropical mathematics. For instance, both addition and
 multiplication are {\em commutative}:
 $$ x \, \oplus \, y \,\,\, = \,\,\, y \,\oplus \, x \qquad \hbox{and}
\qquad
      x \, \odot \, y \,\,\, = \,\,\, y \,\odot \,x. $$
  The {\em distributive law} holds for tropical addition and
  tropical multiplication:
    $$ x \, \odot \, ( y \, \oplus\, z) \quad = \quad
        x \,\odot \, y \,\,\, \oplus \,\,\, x \odot \, z . $$
  Here is a numerical example to show distributivity:
\begin{eqnarray*}
  & 3 \, \odot \, (7 \,\oplus \,11) \,\quad = \quad 3 \odot 7 \,\, \,\,=
\,\,\, 10 ,\\
  & 3 \, \odot 7 \,\,\, \oplus \,\,\, 3 \,\odot 11
  \,\,\, = \,\,\, 10 \,\oplus \, 14 \,\,\, = \,\,\, 10 .
 \end{eqnarray*}
 Both arithmetic operations have a neutral element.
Infinity is the {\em neutral element} for addition and
zero is the {\em neutral element} for multiplication:
$$ x \,\, \oplus \,\,\infty \,\, = \,\, x \qquad \hbox{and} \qquad
     x \,\, \odot \,\, 0  \,\, = \,\, 0. $$
Elementary school students tend to prefer tropical arithmetic
because the multiplication table is easier to memorize,
and even long division becomes easy.
Here are the tropical {\em addition table}
and the tropical {\em  multiplication table}:
$$ \begin{matrix}
\oplus & {\bf 1} & {\bf 2} & {\bf 3} & {\bf 4} & {\bf 5} & {\bf 6} & {\bf
7}          \\
 {\bf 1} &       1   &   1       &     1      &    1      &    1      &
1     &    1                    \\
  {\bf 2} &      1   &    2      &    2      &    2       &    2      &
2    &    2                   \\
  {\bf 3} &      1  &     2      &    3       &    3      &     3     &
3     &    3                 \\
  {\bf 4} &      1 &     2       &   3       &     4      &     4    &
4     &     4                 \\
  {\bf 5} &      1   &    2      &    3       &    4      &     5     &
5     &     5                 \\
  {\bf 6} &       1  &    2       &    3     &      4    &       5    &
6   &      6                \\
  {\bf 7} &      1   &     2      &    3      &      4  &      5    &
6  &        7               \\
 \end{matrix}
 \qquad \qquad  \qquad
 \begin{matrix}
\odot & {\bf 1} & {\bf 2} & {\bf 3} & {\bf 4} & {\bf 5} & {\bf 6} & {\bf
7}          \\
 {\bf 1} &     2   &       3    &      4     &     5     &    6       &
7      &      8                \\
  {\bf 2} &     3    &    4       &    5       &    6     &   7        &
8       &      9                \\
  {\bf 3} &      4   &    5       &    6      &     7    &     8    &
9     &      10                \\
  {\bf 4} &      5   &     6     &     7     &       8    &     9    &
10    &       11               \\
  {\bf 5} &      6  &      7     &     8      &      9    &   10    &   11
&      12                \\
  {\bf 6} &      7  &      8     &    9       &      10   &   11    &   12
&       13              \\
  {\bf 7} &       8  &      9     &    10     &    11   &     12     & 13
&     14             \\
 \end{matrix}
$$
But watch out: tropical arithmetic is tricky
when it comes to subtraction. There is no
 $x$ which we can call
``$ 10$ minus $3$'' because the equation
$\, 3 \oplus x \, = \, 10 \,$ has no solution $x$ at all.
To stay on safe ground,  in this lecture, we shall content ourselves
with using addition $\oplus$ and multiplication $\odot$ only.

It is extremely important to remember that ``$0''$ is the multiplicatively
neutral element. For instance, the tropical {\em Pascal's triangle} looks
like this:
$$ \begin{matrix} &   &  && 0 &&   & &      \\
                         &     &    &0 && 0&    & &    \\
                         &      &    0 && 0 && 0  & & \\
                       &    0 && 0 && 0 && 0 & \\
                        0 && 0 && 0 && 0 && 0 \\
                 \cdots  &  \cdots && \cdots && \cdots && \cdots & \cdots
\\
                \end{matrix}
$$
The rows of Pascal's triangle are the coefficients appearing in
the {\em Binomial Theorem}. For instance, the third row in the triangle
represents the identity
\begin{eqnarray*} (x \oplus y)^3 \quad = &
 (x \oplus y)
\odot
(x \oplus y)
\odot
(x \oplus y) \\
= & \quad
0 \odot x^3 \,\,\oplus \,\, 0 \odot x^2 y
\,\,\oplus \,\, 0 \odot x y^2 \,\, \oplus\, \, 0 \odot y^3.
\end{eqnarray*}
Of course, the zero coefficients can be dropped in this identity:
$$ (x \oplus y)^3 \quad = \quad  x^3 \,\,\oplus \,\,  x^2 y
\,\,\oplus \,\, x y^2 \,\, \oplus\, \, y^3. $$
Moreover, the {\em Freshman's Dream} holds for all powers
in tropical arithmetic:
$$ (x \oplus y)^3 \quad = \quad  x^3 \,\,\,  \oplus\, \,\, y^3. $$
The validity of the three displayed identities is easily verified
by noting that the following equations hold  in classical arithmetic
for all $x,y \in \rr$:
$$ 3 \cdot {\rm min} \{x,y\} \quad = \quad
{\rm min}\{ 3x,2x+y,x+2y,3y \} \quad = \quad
{\rm min}\{ 3x, 3y \}. $$
\smallskip

\noindent {\bf Research problem}: The set of convex polyhedra in
$\rr^n$ can be made into a semiring by taking
$\odot$ as ``Minkowski sum'' and $\oplus$ as
``convex hull of the union''. A natural subalgebra
is the set of all polyhedra which have a fixed
{\em recession cone} $C$. If $n = 1$
and $C = \rr_{\geq 0}$ then we get the
tropical semiring. Develop linear algebra
and algebraic geometry over these semirings,
and implement
efficient software for doing arithmetic
with polyhedra when $n \leq 4$.

 \section{Polynomials}

Let $x_1, x_2 , \ldots , x_n$ be variables which represent
elements in the tropical semiring
$\,(\rr \cup \{\infty\}, \oplus, \odot)$.
A {\em monomial} is any product of these
variables, where repetition is allowed. (Technical note: we allow negative
integer exponents.) By commutativity, we can sort the product and write
monomials in the usual notation, with the variables
raised to exponents:
$$  x_2 \odot x_1 \odot x_3 \odot x_1 \odot x_4 \odot x_2 \odot x_3 \odot
x_2
\quad = \quad x_1^2 x_2^3 x_3^2 x_4. $$
A monomial represents a function from $\rr^n $ to $\rr$.
When evaluating this function in classical arithmetic, what we get
is a linear function:
$$   x_2 +  x_1 + x_3 + x_1 + x_4 + x_2 + x_3 + x_2
\quad = \quad 2 x_1 + 3 x_2 + 2  x_3 + x_4. $$
Every linear function with integer coefficients arises in this manner.

\begin{fact}
Tropical monomials are the linear functions with integer coefficients.
\end{fact}

A {\em tropical polynomial} is a finite
linear combination of tropical monomials:
$$ p(x_1,\ldots,x_n) \quad = \quad
a \odot x_1^{i_1} x_2^{i_2} \cdots x_n^{i_n} \,\,\oplus \,\,
b \odot x_1^{j_1} x_2^{j_2} \cdots x_n^{j_n} \,\,\oplus \,\, \cdots
$$
Here the coefficients $a,b,\ldots$ are real numbers and
the exponents $i_1,j_1,\ldots$ are integers.
Every tropical polynomial represents a function $\,\rr^n \rightarrow \rr$.
When evaluating this function in classical arithmetic, what we get
is the minimum of a finite collection of linear functions, namely,
$$ p(x_1,\ldots,x_n) \quad = \quad
{\rm min} \bigl(
a + i_1 x_1  + \cdots + i_n x_n \,, \,
b + j_1 x_1  + \cdots + j_n x_n \,,\, \ldots \,\bigr) $$
This function $\,p : \rr^n \rightarrow \rr \, $  has
the following three important properties:
\begin{itemize}
\item $p$ is continuous,
\item $p$ is piecewise-linear, where the number of pieces is finite, and
\item $p$ is concave, i.e., $p(\frac{x+y}{2}) \geq \frac{1}{2} (p(x) +
p(y))$ for all $x,y \in \rr^n$.
\end{itemize}
It is known that every function which satisfies these three properties
can be represented as the minimum of a finite set of linear functions. We
conclude:

\begin{fact}
The tropical polynomials in $n$ variables $x_1,\ldots,x_n$
are precisely the piecewise-linear concave functions on $\rr^n$
with integer coefficients.
\end{fact}

As a first example consider the general cubic polynomial in one
variable~$x$,
\begin{equation}
\label{cubic} p(x) \quad = \quad a \odot x^3 \,\,\oplus \,\, b \odot x^2
\,\, \oplus \,\,
c \odot x \,\,\oplus \,\, d .
\end{equation}
To graph this function we draw four lines in the $(x,y)$ plane:
$\, y = 3x + a $,  $\, y =   2x + b $,
$\, y = x + c \,$ and the horizontal line $y= d$.
The value of $p(x)$ is the smallest $y$-value such that
$(x,y)$ is on one of these four lines, i.e.,
the graph of $p(x)$ is the lower envelope of the
lines. All four lines actually contribute if
\begin{equation}
\label{rootsofcubic} b-a \,\, \leq \,\, c-b \,\, \leq d-c .
\end{equation}
These three values of $x$ are the breakpoints where  $p(x)$
fails to be linear, and the cubic  has a corresponding
factorization into three linear factors:
\begin{equation}
\label{factoredcubic}
p(x) \quad = \quad a \odot
 (x \,\oplus \, (b-a)) \odot   (x \,\oplus \,(c-b))  \odot  (x \,\oplus
\,(d-c)) .
 \end{equation}
See Figure 1 for the graph and the roots of the cubic polynomial $p(x)$.

\begin{figure}
\centerline{\scalebox{0.7}{\includegraphics{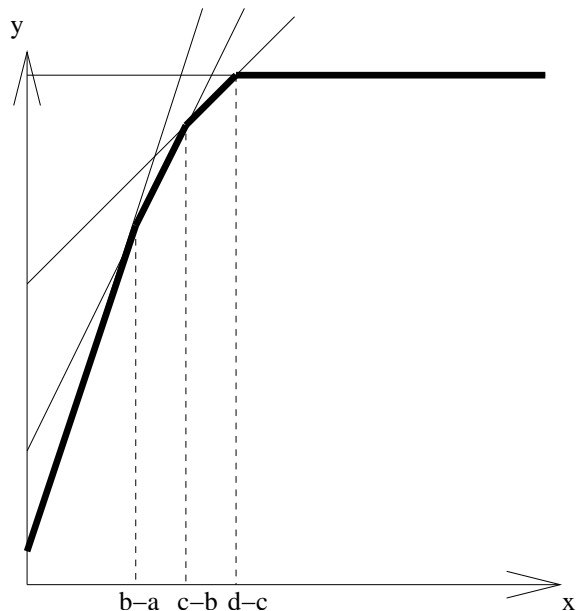}}}
\caption{The graph of a cubic polynomial and its roots}
\end{figure}

Every tropical polynomial function can be written uniquely as a tropical
product of
tropical linear functions (i.e., the {\em
Fundamental Theorem of Algebra} holds tropically). In this statement
we must underline the word ``function''.  Distinct polynomials
can represent the same function. We are not claiming that
every polynomial factors into linear functions. What we are
claiming is that every polynomial can be replaced
by an equivalent polynomial, representing the
same function, that can be factored into linear factors, e.g.,
$$ x^2 \,\,\oplus \,\, 17 \odot x \,\, \oplus \,\, 2 \quad = \quad
     x^2 \,\,\oplus \,\, 1 \odot x \,\, \oplus \,\, 2 \quad = \quad
     (x \, \oplus \, 1)^2 .
$$
Unique factorization of polynomials no longer holds
in two or more variables. Here the situation is more
interesting. Understanding it is our next problem.

\bigskip

\noindent {\bf Research problem}: The factorization of multivariate
tropical polynomials
into irreducible tropical polynomials is not unique.
Here is a simple example:
\begin{eqnarray*}
&
(0 \odot x \,\oplus \, 0) \,\, \odot \,\,
(0 \odot y \,\oplus\, 0)\,\, \odot \,\,
(0 \odot x \odot  y  \,\oplus\, 0)  \\
= &
(0 \odot x \odot  y  \,\oplus\, 0 \odot x \, \oplus \, 0) \,\odot \,
(0 \odot x \odot  y  \,\oplus\, 0 \odot y \, \oplus \, 0) .
\end{eqnarray*}
Develop an algorithm
(with implementation and complexity analysis) for computing
all the irreducible factorizations of
a given tropical polynomial.
Gao and Lauder \cite{GL} have shown
the importance of tropical factorization for the
problem of factoring multivariate polynomials in the
classical sense.

\section{Curves}

A tropical polynomial function $p : \rr^n \rightarrow \rr$
is given as  the minimum of a finite set of linear functions.
We define the {\em hypersurface} $\mathcal{H}(p)$
to be the set of all points $ x \in \rr^n$
at which this minimum is attained at least twice.
Equivalently, a point $x \in \rr^n$ lies in
 $\mathcal{H}(p)$ if and only if $p$ is not linear at $x$.
 For example, if $n = 1$ and $p$ is the cubic in
(\ref{cubic}) with the assumption (\ref{rootsofcubic}), then
$$ \mathcal{H}(p) \quad = \quad
\bigl\{ \, b-a ,\,  c-b , \, d-c \, \bigr\}. $$
Thus the hypersurface
$\mathcal{H}(p)$ is the set of ``roots'' of the polynomial $p(x)$.

 In this section we consider the case of a
 polynomial in two variables:
 $$ p(x,y) \quad = \quad
 \bigoplus_{(i,j)} c_{ij} \odot x^i \odot y^j . $$

\begin{fact}
The tropical curve $\mathcal{H}(p)$ is a finite
graph which is embedded in the plane $\rr^2$.
It has both bounded and unbounded edges,
all edge directions are rational,  and this graph
satisfies a {\em zero tension condition} around each node.
\end{fact}

 The zero tension condition is the following geometric condition.
 Consider any node $(x,y)$ of the graph and suppose it is
 the origin, i.e., $(x,y) = (0,0)$. Then the edges adjacent to this node
 lie on lines with rational slopes. On each such ray 
 emanating from the origin consider the first non-zero
 lattice vector. {\em Zero tension} at $(x,y)$  means that
 the sum of these  vectors is zero.

 Our first example is a {\em line} in the plane. It is defined by a
polynomial:
 $$ p (x,y) \quad = \quad a \odot x \,\,\oplus \,\, b \odot y \,\, \oplus
\,\, c
 \qquad \quad \hbox{where} \,\, a,b,c \in \rr. $$
 The curve $\mathcal{H}(p)$ consists of all
 points $(x,y)$ where the function
 $$ p \,\, : \,\, \rr^2 \,\rightarrow \, \rr \, , \qquad
 (x,y) \,\, \mapsto \,\,{\rm min} \,\bigl( \, a+x ,\, b+y, \, c \,\bigr)
$$
 is not linear. It consists of three
  half-rays emanating from the point $(x,y) = (c-a,c-b)\,$ into
northern, eastern and southwestern direction.

 Here is  a general method for drawing a tropical curve
 $\mathcal{H}(p)$ in the plane. Consider any term
  $\, \gamma \odot x^i \odot y^j \,$ appearing in the polynomial $p$.
  We represent this term by the point $\,(\gamma,i,j)\,$ in $\rr^3$,
  and we compute the convex hull of these points in $\rr^3$.
  Now project the lower envelope of that convex hull into
  the plane under the map $\,\rr^3 \rightarrow \rr^2,\,
  (\gamma,i,j) \mapsto (i,j)$. The image is a planar
convex  polygon together with a distinguished subdivision $\Delta$
into smaller polygons. The tropical curve $\mathcal{H}(p)$
is the dual graph to this subdivision.

\begin{figure}
\centerline{\includegraphics{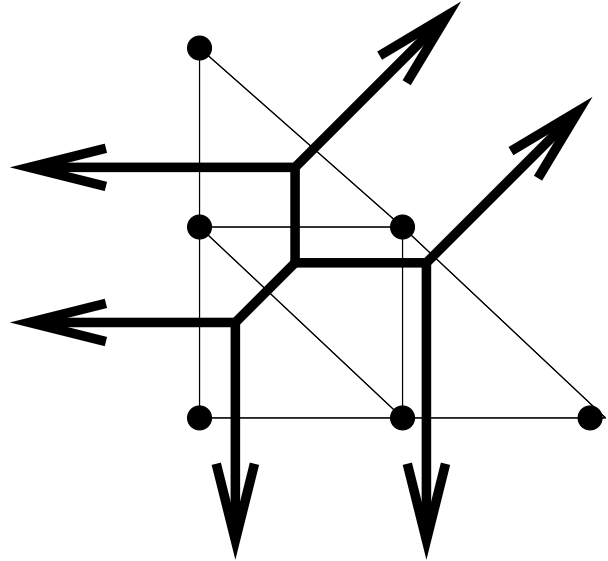}}
\caption{The subdivision $\Delta$ and the tropical curve}
\end{figure}

As an example we consider the general quadratic polynomial
 $$ p(x,y) \quad = \quad a \odot x^2 \,\,\oplus \,\, b \odot xy \,\,
\oplus \,\, c \odot y^2
 \,\, \oplus \,\, d \odot x \,\,\oplus \,\, e \odot y \,\,\oplus \,\, f
.$$
 Then $\Delta$ is a subdivision
 of the triangle with vertices $(0,0)$, $(0,2)$ and $(2,0)$.
 The lattice points $(0,1)$, $(1,0)$, $(1,1)$ are allowed
 to be used as vertices in these subdivisions. Assuming that
 $a,b,c,d,e,f \in \rr$ are general solutions of
 $$ 2b \leq a+c \,  ,\,\,  2d \leq a+f \, ,\,\, 2e \leq c+f, $$
 the subdivision $\Delta$ consists of four triangles,
 three interior edges and six boundary edges.
 The curve $\mathcal{H}(p)$ has four vertices,
 three bounded edges and six half-rays
 (two northern, two eastern and two southwestern). In Figure 2,
$\mathcal{H}(p)$ is shown in bold, and the subdivision $\Delta$ is shown
in
thin lines.

\begin{fact} Tropical curves intersect and interpolate like algebraic
curves do.
\begin{enumerate}
\item Two general lines meet in one point, a line and a
quadric meet in two points, two quadrics meet in
four points, etc....
\item Two general points lie on a unique line,
five general points lie on a unique quadric, etc...
\end{enumerate}
\end{fact}

For a general discussion of {\em B\'ezout's Theorem} in
tropical algebraic geometry, and for many pictures illustrating
Fact 4, we refer to the article \cite{RGST}.

\bigskip

\noindent {\bf Research problem}:
Classify all combinatorial types of
{\em tropical curves in $3$-space} of degree $d$.
Such a curve is a  finite embedded graph
of the form
$$ C \quad = \quad
\mathcal{H}(p_1) \,\cap \,
\mathcal{H}(p_2) \,\cap \,\,\cdots \,\,\cap
\mathcal{H}(p_r) \quad \subset \quad \rr^3 ,$$
where the $p_i$ are tropical polynomials,
 $C$ has $d$ unbounded parallel
halfrays in each of the four coordinate directions, and
all other edges of $C$ are bounded.

\section{Phylogenetics}

An important problem in computational biology is to
construct a {\em phylogenetic tree} from distance
data involving $n$ taxa. These taxa might be organisms
or genes, each represented by a DNA sequence.
For an introduction to phylogenetics we recommend
\cite{Fel} and \cite{SeSt}.
Here is an example, for $n=4$, to illustrate how
such data might arise. Consider an alignment
of four genomes:

$$\begin{matrix}
\hbox{Human:}   & ACAATGTCATTAGCGAT \ldots \\
\hbox{Mouse:}   & ACGTTGTCAATAGAGAT \ldots \\
\hbox{Rat:}     & ACGTAGTCATTACACAT \ldots \\
\hbox{Chicken:} & GCACAGTCAGTAGAGCT \ldots \\
\end{matrix}$$
{}From such sequence data, computational biologists infer the distance
between any
two taxa. There are various algorithms for carrying out this inference.
They are
based on statistical models of evolution. For our discussion, we may
think of the distance between any two strings as a refined version of
the Hamming distance (= the proportion of characters where they differ).
In our (Human, Mouse, Rat, Chicken) example, the inferred distance matrix
might
be the following symmetric $4 \times 4$-matrix:
$$\begin{matrix}
 & H & M & R & C \\
H & 0 & 1.1 & 1.0 & 1.4 \\
M & 1.1 & 0 & 0.3 & 1.3 \\
R & 1.0 & 0.2 & 0 & 1.2 \\
C & 1.4 & 1.3 & 1.2 & 0
\end{matrix}$$
The problem of phylogenetics is to construct a tree with
edge lengths which represent this distance matrix,
provided such a tree exists.
In our example, a tree does exist. It is depicted in Figure 3.
The  number next to the each edge is its length. The distance
between two leaves in the tree is the sum of the lengths
of the edges on the unique path between the two leaves.
For instance, the distance in this tree between ``Human'' and
``Mouse''  equals $\,0.6 + 0.3 + 0.2 \, = \, 1.1$, which is the
corresponding entry in the $4 \times 4$-matrix.

\begin{figure}\label{Tree}
\centerline{\includegraphics{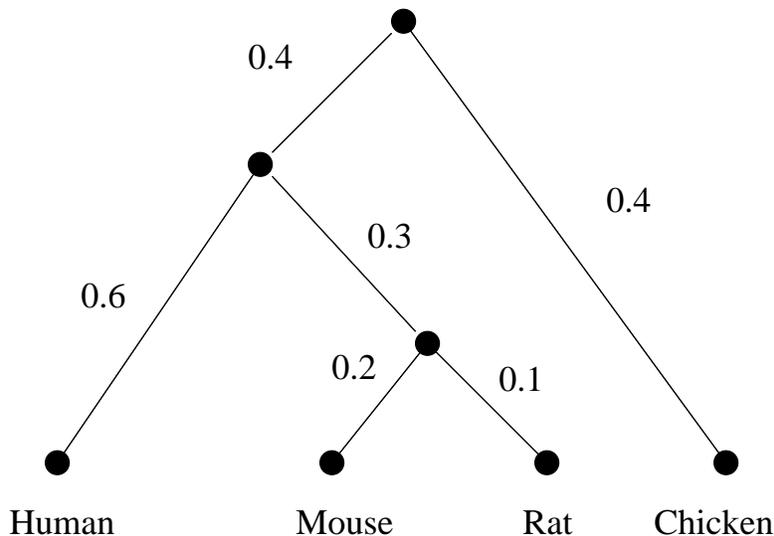}}
\caption{A Phylogenetic Tree}
\end{figure}

In general, considering $n$ taxa,
the {\em distance}
between taxon $i$ and taxon $j$ is  a positive
real number $d_{ij}$ which has been determined
by some bio-statistical method. So, what  we are given is
a  real symmetric $n \times n$-matrix
$$ D \quad = \quad
\begin{pmatrix}
     0     & d_{12}  & d_{13} & \cdots & d_{1n} \\
d_{12} &    0        & d_{23} & \cdots & d_{2n} \\
d_{13} & d_{23} &   0         & \cdots & d_{3n} \\
   \vdots &  \vdots &    \vdots         & \ddots & \vdots \\
 d_{1n} & d_{2n} & d_{3n} & \cdots & 0
\end{pmatrix}.
$$

We may assume that $D$ is a {\em metric},
i.e., the triangle inequalities
$\, d_{ik} \,\leq \ d_{ij} + d_{jk} \,$
hold for all $i,j,k $.
This can be expressed by matrix multiplication:

\begin{fact}
The matrix $D$ represents a metric if and only if
$\,D \odot D = D$.
\end{fact}

We say that $D$ is a {\em tree metric} if there exists
a tree $T$ with $n$ leaves, labeled $1,2,\ldots,n$,
and a positive length for each edge of $T$, such that the
distance from leaf $i$ to leaf $j$ equals $d_{ij}$
for all $i,j$. Tree metrics occur naturally in biology because
they model an evolutionary process  that led to the $n$ taxa.

 Most metrics $D$ are not tree metrics.
 If we are given a metric $D$ that arises
 from some biological data then it is reasonable
 to assume that there exists a tree metric $D_T$
 which is close to $D$. Biologists use a  variety of algorithms
(e.g. ``neighbor joining'') to construct
such a nearby tree $T$ from the given data $D$.
In what follows we state a tropical
characterization of tree metrics.

Let $X = (X_{ij})$ be a symmetric matrix with zeros on the diagonal
whose ${n \choose 2}$ distinct off-diagonal entries are unknowns.
For each quadruple $\{i,j,k,l\}$ $ \subset \{1,2,\ldots,n\}$ we consider
the
following tropical polynomial of degree two:
\begin{equation}
\label{plucker}
p_{i j k l} \quad = \quad
X_{i j} \odot X_{k l} \,\, \oplus \,\,
X_{i k} \odot X_{j l} \,\, \oplus \,\,
X_{i l} \odot X_{j k } .
\end{equation}
This polynomial is the {\em tropical Grassmann-Pl\"ucker relation}.
It defines a hypersurface $\mathcal{H}(p_{i j k l})$ in the
space $\rr^{n \choose 2}$. The {\em tropical Grassmannian}
is the intersection of these ${n \choose 4}$ hypersurfaces.
This is a polyhedral fan, denoted
$$ Gr_{2,n} \quad = \quad \bigcap_{1 \leq i < j<k<l \leq n}
\!\! \mathcal{H}(p_{i j k l} ) \qquad \subset \quad \rr^{n \choose 2}. $$

\begin{fact} A metric $D$ on $ \{1,2,\ldots,n\}$ is a tree metric if and
only if
its negative $X = - D$ is a point in the tropical Grassmannian $Gr_{2,n}$.
\end{fact}

The statement is a reformulation of the
{\em Four Point Condition} in phylogenetics, which states that
$D$ is a tree metric if and only if, for all $1 \leq i < j<k<l \leq n$,
the {\em maximum} of the three numbers
$\,D_{ij} + D_{kl}$, $\,D_{ik} + D_{jl}$ and
$\,D_{il} + D_{jk}$ is attained at least twice.
For $X = -D$, this means that the {\em minimum} of
the three numbers $\,X_{ij} + X_{kl}$, $\,X_{ik} + X_{jl}$ and
$\,X_{il} + X_{jk}$ is attained at least twice,
or, equivalently, $\, X \in \mathcal{H}(p_{i j k l} ) $.
The tropical Grassmannian $Gr_{2,n}$ is also known as the
{\em space of phylogenetic trees} \cite{BHV}. The combinatorial structure
of this
beautiful space is well-studied and well-understood.

Our research suggestion for this section concerns a certain reembedding
of the tropical Grassmannian $Gr_{2,n}$ into a
higher-dimensional space.

\bigskip

\noindent {\bf Research problem}:
Let $n \geq 5$ and consider a metric $D$ on
$\{1,2,\ldots,n\}$.
The {\em triple weights} of the metric $D$
are defined as follows:
$$  D_{ijk} \quad := \quad
D_{ij} + D_{ik} + D_{jk}
\qquad \quad (1 \leq i < j < k \leq n ). $$
This formula specifies a linear map
$\,\psi \,:\, \rr^{n \choose 2} \, \rightarrow \, \rr^{n \choose 3}$.
Our problem is to characterize the image $\, \psi(Gr_{2,n})\,$
of tree space $Gr_{2,n}$ under this linear map.

This problem arises from phylogenies of
sequence alignments in genomics. For such taxa,
it can be more reliable statistically to estimate
the triple weights $D_{ijk}$
rather than the pairwise distances $D_{ij}$.
Pachter and Speyer \cite{PSpey}  showed
that tree metrics are uniquely determined by
their triple weights, i.e., the map
from $\,Gr_{2,n}\,$ onto  $\, \psi(Gr_{2,n})\,$
is a bijection. Find a natural system of
tropical polynomials which define
$\, \psi(Gr_{2,n})\,$ as a tropical subvariety of $\rr^{n \choose 3}$.

\section{Linear Spaces}

This section is concerned with tropical linear spaces.
Generalizing the notion of a line from Section 3,
we define a \emph{tropical hyperplane} to be a subset of
$\rr^n$ which has the form $\mathcal{H}(\ell)$,
where $\ell$ is a tropical linear form in $n$ unknowns:
$$ \ell(x) \quad = \quad a_1 \odot x_n \,\, \oplus \,\,
a_2 \odot x_2 \,\,\oplus \,\, \cdots \,\, \oplus \,\, a_n \, \odot \, x_n
. $$
Here  $a_1$, \dots $a_n$ are arbitrary real constants.
Solving linear equations in tropical mathematics means
computing the intersection of finitely many
hyperplanes $\mathcal{H}(\ell)$. It is tempting to
define tropical linear spaces simply as intersections
of tropical hyperplanes. However, this would not be
a good definition because such arbitrary intersections are not always
pure dimensional, and they do not behave the way
 linear spaces do in classical geometry.

A better notion of tropical linear space is derived
by allowing only those intersections of hyperplanes
which are ``sufficiently complete''. In what follows
we offer a definition which is a direct generalization
of our discussion in Section 4. The idea is that
phylogenetic trees are lines in tropical
projective space, whose Pl\"ucker coordinates
$X_{ij}$ are the negated pairwise distances $d_{ij}$.

We consider  the ${n \choose d}$-dimensional
 space $ \rr ^{n \choose d} \,$ whose
coordinates $X_{i_1 \cdots i_d }$ are indexed
by $d$-element subsets $\{i_1, \ldots,i_d\}$
of $\{1,2,\ldots,n\}$. Let $S$ be any $(d-2)$-element
subset of $\{1,2,\ldots,n\}$ and  let $i$, $j$, $k$ and $l$ be any four
distinct
indices in $\{ 1, \ldots, n \} \backslash S$.
The corresponding {\em three-term Grassmann Pl\"ucker relation}
$\,p_{S,  ij k l}\, $ is the
following tropical polynomial of degree two:
\begin{equation}
\label{plucker2}
p_{S,  ij k l} \quad = \quad
X_{S i j} \odot X_{S k l} \,\,\, \oplus \,\,\,
X_{S i k} \odot X_{S j l} \,\,\, \oplus \,\,\,
X_{S i l} \odot X_{S j k } .
\end{equation}
We define the \emph{three-term tropical Grassmannian} to be the
intersection
$$ Gr_{d.n} \,\,\,= \,
\bigcap_{S, i,j,k,l} \mathcal{H}(p_{S,ijkl}) \,\quad \subset \quad \,
\rr^{\binom{n}{d}},$$
where the intersection is over all $S$, $i$, $j$, $k$, $l$ as above.
Note that in the special case $d=2$ we have $S = \emptyset$,
the polynomial (\ref{plucker2}) is the
four point condition in (\ref{plucker}),
and $Gr_{d,n}$ is the space of trees
which was discussed in Section 4.

We now fix an arbitrary point $\, X = (X_{i_1 \cdots i_d})\,$ in the
three-term tropical Grassmannian $\, Gr_{d.n} $.
For any $(d+1)$-subset $\{j_0,j_1,\ldots,j_{d}\}$ of
$\{1,2,\ldots,n\}$ we consider the following tropical linear form
in the  variables $x_1,\ldots,x_n$:
\begin{equation}
\label{cramer} \ell^X_{j_0 j_1 \cdots j_d} \quad = \quad
\bigoplus_{r=0}^{d} X_{j_0 \cdots \widehat{j_r} \cdots j_{d}} \odot x_r.
\end{equation}
The $\widehat{\phantom{j_r}}$ means to omit $j_r$. The {\em tropical
linear space} associated with the point $X $
 is the following set:
$$ L_X \quad = \quad \bigcap \mathcal{H}( \ell^X_{j_0 j_1 \cdots j_n} )
\qquad \subset \quad \rr^n . $$
Here the intersection is over all $(d+1)$-subsets
$\{j_0,j_1,\ldots,j_{d}\}$ of
$\{1,2,\ldots,n\}$.

The tropical linear spaces are precisely the sets $L_X$ where
$X$ is any point in $Gr_{d,n} \subset \rr ^{n \choose d}$.
The ``sufficient completeness'' referred to in the
first paragraph of this section means that we need to
solve linear equations using {\em Cramer's rule}, in all
possible ways, in order for an intersection of
hyperplanes to actually be a linear space.
The definition of linear space given here is
 more inclusive than the one used in
\cite{DSS, RGST, SS}, where $L_X$ was required to
come from ordinary algebraic geometry over  a field
with a suitable valuation.

For example, a $3$-dimensional tropical linear subspace of
$\rr^n$ (a.k.a.~two-dimensional plane in
tropical projective $(n-1)$-space)
is the intersection of ${n \choose 4}$ tropical
hyperplanes, each of whose defining linear form
has four terms:
$$ \ell^X_{j_0 j_1 j_2 j_3} \,\,\, = \,\,\,\,
X_{j_0 j_1 j_2} \odot x_{j_3} \,\, \oplus \,\,
X_{j_0 j_1 j_3} \odot x_{j_2} \,\, \oplus \,\,
X_{j_0 j_2 j_3} \odot x_{j_1} \,\, \oplus \,\,
X_{j_1 j_2 j_3} \odot x_{j_0} . $$

We note that even the very special case when each coordinate of $X$
is either $0$ (the multiplicative unit) or $\infty$ (the additive unit)
is really interesting. Here $L_X$ is a polyhedral fan known as the
 {\em Bergman fan} of a matroid \cite{AK, Stu}.

Tropical linear spaces have many of the properties of ordinary linear
spaces.
First, they are pure polyhedral complexes of the correct dimension:

\begin{fact}
Each maximal cell of
the tropical linear space $L_X$ is $d$-dimensional.
\end{fact}

Every tropical linear space $L_X$ determines its vector of tropical
Pl\"ucker coordinates
$X$ uniquely up to tropical multiplication (= classical addition) by a
common scalar. If $L$ and $L'$ are tropical linear spaces of dimensions
$d$ and $d'$ with $d+d' \geq n$, then $L$ and $L'$ meet. It is not quite
true that two tropical linear spaces intersect in a tropical linear space
but it is almost true. If $L$ and $L'$ are tropical linear spaces of
dimensions $d$ and $d'$ with $d+ d' \geq n$ and $v$ is a generic small
vector then $L \cap (L' + v)$ is a tropical linear space of dimension
$d+d' -n$. Following \cite{RGST}, it makes sense to define the {\em stable
intersection} of $L$ and $L'$ by taking the limit of $L \cap (L'+v)$
as $v$ goes to zero, and this limit will again be a tropical linear space
of dimension $d+d'-n$.

\bigskip

\noindent \textbf{Research Problem}:
It is not true that a $d$-dimensional tropical linear space can always be
written as the intersection of $n-d$ tropical hyperplanes. The definition
shows that $\binom{n}{d +1 }$
hyperplanes are always enough. What is the minimum number of tropical
hyperplanes needed
 to cut out any tropical linear space of dimension $d$ in $n$-space? Are
$n$ hyperplanes always enough?


\thebibliography{99}

\raggedright


 \bibitem{Ard} F.~Ardila: A tropical morphism related to the hyperplane
 arrangement of the complete bipartite graph,
 \texttt{http://www.arXiv.org/math.CO/0404287}.

\bibitem{AK} F.~Ardila and C.~Klivans: The Bergman complex of a matroid
and phylogenetic trees, \texttt{http://www.arXiv.org/math.CO/0311370}.

\bibitem{Berg} G.~Bergman: The logarithmic limit-set of an
algebraic variety, \emph{Trans. of the
AMS} {\bf 157} (1971) 459--469.

\bibitem{BG} R.~Bieri and J.R.J.~Groves:
The geometry of the set of characters induced by valuations.
\emph{J.~reine und angewandte Mathematik}
 {\bf 347} (1984) 168-195.

\bibitem{BHV} L.~Billera, S.~Holmes and K.~Vogtman: Geometry of
the space of phylogenetic trees, \emph{Advances in Applied Mathematics}
{\bf 27} (2001) 733--767.

\bibitem{But} P. Butkovi\v c: Max-algebra: the linear algebra of
combinatorics?  \emph{Linear Algebra Appl.}  \textbf{367}  (2003)
313--335.

\bibitem{Dev1} M.~Develin:  Tropical secant varieties of linear spaces,
{\tt http://www.arXiv.org/math.CO/0405115}.

\bibitem{Dev2} M.~Develin:
The space of $n$ points on a tropical line in $d$-space,
{\tt http://www.arXiv.org/math.CO/0401224}.

\bibitem{DS} M.~Develin and B.~Sturmfels: Tropical convexity,
{\sl Documenta Mathematica} {\bf 9} (2004) 1--27.

\bibitem{DSS} M.~Develin, F.~Santos and B.~Sturmfels: On the rank of a
tropical matrix, {\tt http://www.arXiv.org/math.CO/0312114}.

\bibitem{Fel} J. Felsenstein: {\sl Inferring Phylogenies}, Sinauer
Associates, Inc.,
Sunderland, 2003.

\bibitem{GL} S.~Gao and A.~Lauder: Decomposition of polytopes
and polynomials, {\sl Discrete and Computational Geometry}
{\bf 26} (2001) 89--104.

\bibitem{IKS}
 I. Itenberg, V. Kharlamov and E. Shustin:
 Welschinger invariant and enumeration of real plane rational curves
 {\tt  http://www.arXiv.org/math.AG/0303378}.

 \bibitem{Jos} M.~Joswig: Tropical halfspaces, {\tt
http://www.arXiv.org/math.CO/0312068}.

\bibitem{Mik} G.~Mikhalkin: Enumerative tropical geometry in $\rr^2$,
\texttt{ http://www.arXiv.org/math.AG/0312530}.

\bibitem{Mik2} G.~Mikhalkin:  Amoebas of algebraic varieties and tropical
geometry,
\texttt{ http://www.arXiv.org/math.AG/0403015}.

\bibitem{PSturm} L.~Pachter and B.~Sturmfels: Tropical geometry of
statistical models,
{\tt http://www.arXiv.org/q-bio.QM/0311009}, to appear in
{\sl Proceedings of the National Academy of Sciences}.

\bibitem{PSpey} L. Pachter and D. Speyer: Reconstructing trees from
subtree weights, to appear in \emph{Applied Mathematics Letters},
\texttt{ http://www.arXiv.org/math.CO/0311156}.

\bibitem{Pin} J.-E.~Pin:
Tropical semirings. \emph{Idempotency} (Bristol, 1994), 50--69,
Publ. Newton Inst., {\bf 11}, Cambridge Univ. Press, Cambridge, 1998.

\bibitem{RGST} J.~Richter-Gebert, B.~Sturmfels, and T.~Theobald:
First steps in tropical geometry,  {\tt math.AG/0306366}, to appear in
\emph{Idempotent Mathematics and  Mathematical  Physics},
Proceedings Vienna 2003, (editors G.L.~Litvinov  and  V.P.~Maslov),
American Mathematical Society, 2004.

\bibitem{SeSt} C.~Semple and M.~Steel:  {\em Phylogenetics},
Oxford University Press, Oxford, 2003.

 \bibitem{Shu}  E. Shustin: Patchworking singular algebraic curves,
 non-Archimedean amoebas and enumerative geometry,
 {\tt   http://www.arXiv.org/math.AG/0211278}.

 \bibitem{Simon} I. Simon: Recognizable sets with multiplicities in the
tropical semiring.  Mathematical foundations of computer science,
(Carlsbad, 1988),  107--120, Lecture Notes in Comput. Sci., \textbf{324},
Springer, Berlin, 1988.

\bibitem{SS} D. Speyer and B. Sturmfels: The tropical Grassmannian,
\emph{Advances in Geometry}, \textbf{4}, Issue 3 (2004), p. 389-411

\bibitem{Stu} B. Sturmfels: \emph{Solving Systems of Polynomial
Equations},
American Mathematical Society, CMBS Series, {\bf 97}, 2002.
\end{document}